# ASYMPTOTIC PROPERTIES OF BRIDGE ESTIMATORS IN SPARSE HIGH-DIMENSIONAL REGRESSION MODELS

By Jian Huang,[1] Joel L. Horowitz[2] and Shuangge Ma

*University of Iowa, Northwestern University and Yale University*

We study the asymptotic properties of bridge estimators in sparse, high-dimensional, linear regression models when the number of covariates may increase to infinity with the sample size. We are particularly interested in the use of bridge estimators to distinguish between covariates whose coefficients are zero and covariates whose coefficients are nonzero. We show that under appropriate conditions, bridge estimators correctly select covariates with nonzero coefficients with probability converging to one and that the estimators of nonzero coefficients have the same asymptotic distribution that they would have if the zero coefficients were known in advance. Thus, bridge estimators have an oracle property in the sense of Fan and Li [*J. Amer. Statist. Assoc.* **96** (2001) 1348–1360] and Fan and Peng [*Ann. Statist.* **32** (2004) 928–961]. In general, the oracle property holds only if the number of covariates is smaller than the sample size. However, under a partial orthogonality condition in which the covariates of the zero coefficients are uncorrelated or weakly correlated with the covariates of nonzero coefficients, we show that marginal bridge estimators can correctly distinguish between covariates with nonzero and zero coefficients with probability converging to one even when the number of covariates is greater than the sample size.

**1. Introduction.** Consider the linear regression model

$$Y_i = \beta_0 + \mathbf{x}_i'\boldsymbol{\beta} + \varepsilon_i, \qquad i = 1, \ldots, n,$$

where $Y_i \in \mathbb{R}$ is a response variable, $\mathbf{x}_i$ is a $p_n \times 1$ covariate vector and the $\varepsilon_i$'s are i.i.d. random error terms. Without loss of generality, we assume that $\beta_0 = 0$. This can be achieved by centering the response and covariates.

Received March 2006; accepted April 2007.
[1]Supported in part by NIH Grant NCI/NIH P30 CA 086862-06.
[2]Supported in part by NSF Grant SES-0352675.

We are interested in estimating the vector of regression coefficients $\boldsymbol{\beta} \in \mathbb{R}^{p_n}$ when $p_n$ may increase with $n$ and $\boldsymbol{\beta}$ is sparse in the sense that many of its elements are zero. We estimate $\boldsymbol{\beta}$ by minimizing the penalized least squares objective function

$$(1) \qquad L_n(\boldsymbol{\beta}) = \sum_{i=1}^{n}(Y_i - \mathbf{x}_i'\boldsymbol{\beta})^2 + \lambda_n \sum_{j=1}^{p_n} |\beta_j|^\gamma,$$

where $\lambda_n$ is a penalty parameter. For any given $\gamma > 0$, the value $\widehat{\boldsymbol{\beta}}_n$ that minimizes (1) is called a bridge estimator [Frank and Friedman (1993) and Fu (1998)]. The bridge estimator includes two important special cases. When $\gamma = 2$, it is the familiar ridge estimator [Hoerl and Kennard (1970)]. When $\gamma = 1$, it is the LASSO estimator [Tibshirani (1996)], which was introduced as a variable selection and shrinkage method. When $0 < \gamma \leq 1$, some components of the estimator minimizing (1) can be exactly zero if $\lambda_n$ is sufficiently large [Knight and Fu (2000)]. Thus, the bridge estimator for $0 < \gamma \leq 1$ provides a way to combine variable selection and parameter estimation in a single step. In this article we provide conditions under which the bridge estimator for $0 < \gamma < 1$ can correctly distinguish between nonzero and zero coefficients in sparse high-dimensional settings. We also give conditions under which the estimator of the nonzero coefficients has the same asymptotic distribution that it would have if the zero coefficients were known with certainty.

Knight and Fu (2000) studied the asymptotic distributions of bridge estimators when the number of covariates is finite. They showed that, for $0 < \gamma \leq 1$, under appropriate regularity conditions, the limiting distributions can have positive probability mass at 0 when the true value of the parameter is zero. Their results provide a theoretical justification for the use of bridge estimators to select variables (i.e., to distinguish between covariates whose coefficients are exactly zero and covariates whose coefficients are nonzero). In addition to bridge estimators, other penalization methods have been proposed for the purpose of simultaneous variable selection and shrinkage estimation. Examples include the SCAD penalty [Fan (1997) and Fan and Li (2001)] and the Elastic-Net (Enet) penalty [Zou and Hastie (2005)]. For the SCAD penalty, Fan and Li (2001) studied asymptotic properties of penalized likelihood methods when the number of parameters is finite. Fan and Peng (2004) considered the same problem when the number of parameters diverges. Under certain regularity conditions, they showed that there exist local maximizers of the penalized likelihood that have an oracle property. Here the oracle property means that the local maximizers can correctly select the nonzero coefficients with probability converging to one and that the estimators of the nonzero coefficients are asymptotically normal with the same means and covariances that they would have if the zero coefficients were known in advance. Therefore, the local maximizers are asymptotically as efficient as the ideal estimator assisted by an oracle who knows which coefficients are nonzero.

Several other studies have investigated the properties of regression estimators when the number of covariates increases to infinity with sample size. See, for example, Huber (1981) and Portnoy (1984, 1985). Portnoy (1984, 1985) provided conditions on the growth rate of $p_n$ that are sufficient for consistency and asymptotic normality of a class of M-estimators of regression parameters. However, Portnoy did not consider penalized regression or selection of variables in sparse models. Bair et al. (2006) proved consistency of supervised principal components analysis under a partial orthogonality condition, but they also did not consider penalized regression. There have been several other studies of large sample properties of high-dimensional problems in settings related to but different from ours. Examples include Van der Laan and Bryan (2001), Bühlmann (2006), Fan, Peng and Huang (2005), Huang, Wang and Zhang (2005), Huang and Zhang (2005) and Kosorok and Ma (2007). Fan and Li (2006) provide a review of statistical challenges in high-dimensional problems that arise in many important applications.

We study the asymptotic properties of bridge estimators with $0 < \gamma < 1$ when the number of covariates $p_n$ may increase to infinity with $n$. We are particularly interested in the use of bridge estimators to distinguish between covariates with zero and nonzero coefficients. Our study extends the results of Knight and Fu (2000) to infinite-dimensional parameter settings. We show that for $0 < \gamma < 1$ the bridge estimators can correctly select covariates with nonzero coefficients and that, under appropriate conditions on the growth rates of $p_n$ and $\lambda_n$, the estimators of nonzero coefficients have the same asymptotic distribution that they would have if the zero coefficients were known in advance. Therefore, bridge estimators have the oracle property of Fan and Li (2001) and Fan and Peng (2004). The permitted rate of growth of $p_n$ depends on the penalty function form specified by $\gamma$. We require that $p_n < n$; that is, the number of covariates must be smaller than the sample size.

The condition that $p_n < n$ is needed for identification and consistent estimation of the regression parameter. While this condition is often satisfied in applications, there are important settings in which it is violated. For example, in studies of relationships between a phenotype and microarray gene expression profiles, the number of genes (covariates) is typically much greater than the sample size, although the number of genes that are actually related to the clinical outcome of interest is generally small. Often a goal of such studies is to find these genes. Without any further assumption on the covariate matrix, the regression parameter is in general not identifiable if $p_n > n$. However, if there is suitable structure in the covariate matrix, it is possible to achieve consistent variable selection and estimation. A special case is when the columns of the covariate matrix $\mathbf{X}$ are mutually orthogonal. Then each regression coefficient can be estimated by univariate regression. But, in practice, mutual orthogonality is often too strong an assumption. Furthermore, when $p_n > n$, mutual orthogonality of all covariates is not possible, since the rank of $\mathbf{X}$ is at most $n-1$. We consider a partial orthogonality condition in which the covariates with zero coefficients are uncorrelated or only weakly correlated with the covariates with nonzero coefficients. We study a univariate version of the bridge estimator. We show that under the partial orthogonality condition and certain other conditions, the marginal bridge estimator can consistently distinguish between zero coefficients and nonzero coefficients even when the number of covariates is greater than $n$, although it does not yield consistent estimation of the parameters. After the covariates with nonzero coefficients are consistently selected, we can use any reasonable method to consistently estimate their coefficients if the number of nonzero coefficients is relatively small, as it is in sparse models. The partial orthogonality condition appears to be reasonable in microarray data analysis, where the genes that are correlated with the phenotype of interest may be in different functional pathways from the genes that are not related to the phenotype [Bair et al. (2006)]. Fan and Lv (2006) also studied univariate screening in high-dimensional regression problems and provided conditions under which it can be used to reduce the exponentially growing dimensionality of a model.

The rest of this paper is organized as follows. In Section 2 we present asymptotic results for bridge estimators with $0 < \gamma < 1$ and $p_n \to \infty$ as $n \to \infty$. We treat a general covariate matrix and a covariate matrix that satisfies our partial orthogonality condition. In Section 3 we present results for marginal bridge estimators under the partial orthogonality condition. In Section 4 simulation studies are used to assess the finite sample performance of bridge estimators. Concluding remarks are given in Section 5. Proofs of the results stated in Sections 2 and 3 are given in Section 6.

**2. Asymptotic properties of bridge estimators.** Let the true parameter value be $\boldsymbol{\beta}_{n0}$. The subscript $n$ indicates that $\boldsymbol{\beta}_{n0}$ depends on $n$, but for simplicity of notation, we will simply write $\boldsymbol{\beta}_0$. Let $\boldsymbol{\beta}_0 = (\boldsymbol{\beta}'_{10}, \boldsymbol{\beta}'_{20})'$, where $\boldsymbol{\beta}_{10}$ is a $k_n \times 1$ vector and $\boldsymbol{\beta}_{20}$ is a $m_n \times 1$ vector. Suppose that $\boldsymbol{\beta}_{10} \neq \mathbf{0}$ and $\boldsymbol{\beta}_{20} = \mathbf{0}$, where $\mathbf{0}$ is the vector with all components zero. So $k_n$ is the number of nonzero coefficients and $m_n$ is the number of zero coefficients. We note that it is unknown to us which coefficients are nonzero and which are zero. We partition $\boldsymbol{\beta}_0$ this way to facilitate the statement of the assumptions.

Let $\mathbf{x}_i = (x_{i1}, \ldots, x_{ip_n})'$ be the $p_n \times 1$ vector of covariates of the $i$th observation, $i = 1, \ldots, n$. We assume that the covariates are fixed. However, we note that for random covariates, the results hold conditionally on the covariates. We assume that the $Y_i$'s are centered and the covariates are standardized, that is,

(2) $$\sum_{i=1}^n Y_i = 0, \qquad \sum_{i=1}^n x_{ij} = 0 \quad \text{and} \quad \frac{1}{n}\sum_{i=1}^n x_{ij}^2 = 1, \qquad j=1,\ldots,p_n.$$

We also write $\mathbf{x}_i = (\mathbf{w}_i', \mathbf{z}_i')'$, where $\mathbf{w}_i$ consists of the first $k_n$ covariates (corresponding to the nonzero coefficients), and $\mathbf{z}_i$ consists of the remaining $m_n$ covariates (those with zero coefficients). Let $\mathbf{X}_n$, $\mathbf{X}_{1n}$ and $\mathbf{X}_{2n}$ be the matrices whose transposes are $\mathbf{X}_n' = (\mathbf{x}_1,\ldots,\mathbf{x}_n)$, $\mathbf{X}_{1n}' = (\mathbf{w}_1,\ldots,\mathbf{w}_n)$ and $\mathbf{X}_{2n}' = (\mathbf{z}_1,\ldots,\mathbf{z}_n)$, respectively. Let

$$\Sigma_n = n^{-1}\mathbf{X}_n'\mathbf{X}_n \quad \text{and} \quad \Sigma_{1n} = n^{-1}\mathbf{X}_{1n}'\mathbf{X}_{1n}.$$

Let $\rho_{1n}$ and $\rho_{2n}$ be the smallest and largest eigenvalues of $\Sigma_n$, and let $\tau_{1n}$ and $\tau_{2n}$ be the smallest and largest eigenvalues of $\Sigma_{1n}$, respectively.

We now state the conditions for consistency and oracle efficiency of bridge estimators with general covariate matrices.

(A1) $\varepsilon_i, \varepsilon_2, \ldots$ are independent and identically distributed random variables with mean zero and variance $\sigma^2$, where $0 < \sigma^2 < \infty$.

(A2) (a) $\rho_{1n} > 0$ for all $n$; (b) $(p_n + \lambda_n k_n)(n\rho_{1n})^{-1} \to 0$.

(A3) (a) $\lambda_n(k_n/n)^{1/2} \to 0$; (b) $\lambda_n n^{-\gamma/2}(\rho_{1n}/\sqrt{p_n})^{2-\gamma} \to \infty$.

(A4) There exist constants $0 < b_0 < b_1 < \infty$ such that

$$b_0 \le \min\{|\beta_{1j}|, 1 \le j \le k_n\} \le \max\{|\beta_{1j}|, 1 \le j \le k_n\} \le b_1.$$

(A5) (a) There exist constants $0 < \tau_1 < \tau_2 < \infty$ such that $\tau_1 \le \tau_{1n} \le \tau_{2n} \le \tau_2$ for all $n$; (b)

$$n^{-1/2} \max_{1 \le i \le n} \mathbf{w}_i'\mathbf{w}_i \to 0.$$

Condition (A1) is standard in linear regression models. Condition (A2)(a) implies that the matrix $\Sigma_n$ is nonsingular for each $n$, but it permits $\rho_{1n} \to 0$ as $n \to \infty$. As we will see in Theorem 2, $\rho_{1n}$ affects the rate of convergence of bridge estimators. Condition (A2)(b) is used in the consistency proof. Condition (A3) is needed in the proofs of the rate of convergence, oracle property and asymptotic normality. To get a better sense of this condition, suppose that $0 < c_1 < \rho_{1n} \le \rho_{2n} < c_2 < \infty$ for some constants $c_1$ and $c_2$ and for all $n$ and that the number of nonzero coefficients is finite. Then (A3) simplifies to

(A3)* (a) $\lambda_n n^{-1/2} \to 0$; (b) $\lambda_n^2 n^{-\gamma} p_n^{-(2-\gamma)} \to \infty$.

Condition (A3)*(a) states that the penalty parameter $\lambda_n$ must always be $o(n^{1/2})$. Suppose that $\lambda_n = n^{(1-\delta)/2}$ for a small $\delta > 0$. Then (A3)*(b) requires that $p_n^{2-\gamma}/n^{1-\delta-\gamma} \to 0$. So the smaller the $\gamma$, the larger $p_n$ is allowed. This condition excludes $\gamma = 1$, which corresponds to the LASSO estimator. If $p_n$ is finite, then this condition is the same as that assumed by Knight and Fu [(2000), page 1361]. Condition (A4) assumes that the nonzero coefficients are uniformly bounded away from zero and infinity. Condition (A5)(a) assumes that the matrix $\Sigma_{1n}$ is strictly positive definite. In sparse problems, $k_n$ is small relative to $n$, so this assumption is reasonable in such problems. Condition (A5)(b) is needed in the proof of asymptotic normality of the estimators of nonzero coefficients. Under condition (A3)(a), this condition is satisfied if all the covariates corresponding to the nonzero coefficients are bounded by a constant $C$. This is because, under (A3)(a), $n^{-1/2} \max_{1 \leq i \leq n} \mathbf{w}_i' \mathbf{w}_i \leq n^{-1/2} k_n C \to 0$.

In the following, the $L_2$ norm of any vector $\mathbf{u} \in R^{p_n}$ is denoted by $\|\mathbf{u}\|$; that is, $\|\mathbf{u}\| = [\sum_{j=1}^{p_n} u_j^2]^{1/2}$.

THEOREM 1 (Consistency). *Let $\widehat{\boldsymbol{\beta}}_n$ denote the minimizer of* (1). *Suppose that $\gamma > 0$ and that conditions* (A1)(a), (A2), (A3)(a) *and* (A4) *hold. Let $h_n = \rho_{1n}^{-1}(p_n/n)^{1/2}$ and $h_n' = [(p_n + \lambda_n k_n)/(n\rho_{1n})]^{1/2}$. Then $\|\widehat{\boldsymbol{\beta}}_n - \boldsymbol{\beta}_0\| = O_p(\min\{h_n, h_n'\})$.*

We note that $\rho_{1n}^{1/2}$ and $\rho_{1n}$ appear in the denominators of $h_{1n}$ and $h_{2n}$, respectively. Therefore, $h_{2n}$ may not converge to zero faster than $h_{1n}$ if $\rho_{1n} \to 0$. If $\rho_{1n} > \rho_1 > 0$ for all $n$, Theorem 1 yields the rate of convergence $O_p(h_{2n}) = O_p((p_n/n)^{1/2})$. If $p_n$ is finite and $\rho_{1n} > \rho_1 > 0$ for all $n$, then the rate of convergence is the familiar $n^{-1/2}$. However, if $\rho_{1n} \to 0$, the rate of convergence will be slower than $n^{-1/2}$.

This result is related to the consistency result of Portnoy (1984). If $\rho_{1n} > \rho_1 > 0$ for all $n$, which Portnoy assumed, then the rate of convergence in Theorem 1 is the same as that in Theorem 3.2 of Portnoy (1984). Here, however, we consider penalized least squares estimators, whereas Portnoy considered general M-estimators in a linear regression model without penalty. In addition, Theorem 1 is concerned with the minimizer of the objective function (1). In comparison, Theorem 3.2 of Portnoy shows that there exists a root of an M-estimating equation with convergence rate $O_p((p_n/n)^{1/2})$.

THEOREM 2 (Oracle property). *Let $\widehat{\boldsymbol{\beta}}_n = (\widehat{\boldsymbol{\beta}}_{1n}, \widehat{\boldsymbol{\beta}}_{2n})$, where $\widehat{\boldsymbol{\beta}}_{1n}$ and $\widehat{\boldsymbol{\beta}}_{2n}$ are estimators of $\boldsymbol{\beta}_{10}$ and $\boldsymbol{\beta}_{20}$, respectively. Suppose that $0 < \gamma < 1$ and that conditions* (A1) *to* (A5) *are satisfied. We have the following:*

(i) $\widehat{\boldsymbol{\beta}}_{2n} = \mathbf{0}$ *with probability converging to* 1.

(ii) *Let $s_n^2 = \sigma^2 \boldsymbol{\alpha}_n' \Sigma_{1n}^{-1} \boldsymbol{\alpha}_n$ for any $k_n \times 1$ vector $\boldsymbol{\alpha}_n$ satisfying $\|\boldsymbol{\alpha}_n\|_2 \leq 1$.
Then*

$$
(3) \quad \begin{aligned} n^{1/2} s_n^{-1} \boldsymbol{\alpha}_n' (\widehat{\boldsymbol{\beta}}_{1n} - \boldsymbol{\beta}_{10}) \\ = n^{-1/2} s_n^{-1} \sum_{i=1}^n \varepsilon_i \boldsymbol{\alpha}_n' \Sigma_{1n}^{-1} \mathbf{w}_i + o_p(1) \to_D N(0,1), \end{aligned}
$$

*where $o_p(1)$ is a term that converges to zero in probability uniformly with respect to $\boldsymbol{\alpha}_n$.*

Theorem 2 states that the estimators of the zero coefficients are exactly zero with high probability when $n$ is large and that the estimators of the nonzero parameters have the same asymptotic distribution that they would have if the zero coefficients were known. This result is stated in a way similar to Theorem 2 of Fan and Peng (2004). Fan and Peng considered maximum penalized likelihood estimation. Their results are concerned with local maximizers of the penalized likelihood. These results do not imply existence of an estimator with the properties of the local maximizer without auxiliary information about the true parameter value that enables one to choose the localization neighborhood. In contrast, our Theorem 2 is for the global minimizer of the penalized least squares objective function, which is a feasible estimator. In addition, Fan and Peng (2004) require that the number of parameters, $p_n$, to satisfy $p_n^5/n \to 0$, which is more restrictive than our assumption for the linear regression model.

Let $\widehat{\beta}_{1nj}$ and $\beta_{10j}$ be the $j$th components of $\widehat{\boldsymbol{\beta}}_{1n}$ and $\boldsymbol{\beta}_{10}$, respectively. Set $\boldsymbol{\alpha}_n = \mathbf{e}_j$ in Theorem 2, where $\mathbf{e}_j$ is the unit vector whose only nonzero element is the $j$th element and let $s_{nj}^2 = \sigma^2 \mathbf{e}_j' \Sigma_{1n}^{-1} \mathbf{e}_j$. Then we have $n^{1/2} s_{nj}^{-1} (\widehat{\beta}_{1nj} - \beta_{10j}) \to_D N(0,1)$. Thus, Theorem 2 provides asymptotic justification for the following steps to compute an approximate standard error of $\widehat{\beta}_{1nj}$: (i) Compute the bridge estimator for a given $\gamma$; (ii) exclude the covariates whose estimates are zero; (iii) compute a consistent estimator $\widehat{\sigma}^2$ of $\sigma^2$ based on the sum of residual squares; (iv) compute $\widehat{s}_{nj}^{-1} = \widehat{\sigma}(\mathbf{e}_j' \Sigma_{1n}^{-1} \mathbf{e}_j)^{1/2}$, which gives an approximate standard error of $\widehat{\beta}_{1nj}$.

Theorem 1 holds for any $\gamma > 0$. However, Theorem 2 assumes that $\gamma$ is strictly less than 1, which excludes the LASSO estimator.

**3. Asymptotic properties of marginal bridge estimators under partial orthogonality condition.** Although the results in Section 2 allow the number of covariates $p_n \to \infty$ as the sample size $n \to \infty$, they require that $p_n < n$. In this section we show that, under a partial orthogonality condition on the covariate matrix, we can consistently identify the covariates with zero coefficients using a marginal bridge objective function, even when the number of covariates increases almost exponentially with $n$. The precise statement of partial orthogonality is given in condition (B2) below. The marginal bridge estimator is computationally simple and can be used to screen out the covariates with zero coefficients, thereby reducing the exponentially growing dimension of the model to a more manageable one. The nonzero coefficients can be estimated in a second step, as is explained later in this section.

The marginal bridge objective function is

$$(4) \qquad U_n(\boldsymbol{\beta}) = \sum_{j=1}^{p_n} \sum_{i=1}^{n} (Y_i - x_{ij}\beta_j)^2 + \lambda_n \sum_{j=1}^{p_n} |\beta_j|^\gamma.$$

Let $\widetilde{\boldsymbol{\beta}}_n$ be the value that minimizes $U_n$. Write $\widetilde{\boldsymbol{\beta}}_n = (\widetilde{\boldsymbol{\beta}}'_{n1}, \widetilde{\boldsymbol{\beta}}'_{n2})'$ according to the partition $\boldsymbol{\beta}_0 = (\boldsymbol{\beta}'_{10}, \boldsymbol{\beta}'_{20})'$. Let $K_n = \{1, \ldots, k_n\}$ and $J_n = \{k_n+1, \ldots, p_n\}$ be the set of indices of nonzero and zero coefficients, respectively. Let

$$(5) \qquad \xi_{nj} = n^{-1} \mathrm{E}\left(\sum_{i=1}^{n} Y_i x_{ij}\right) = n^{-1} \sum_{i=1}^{n} (\mathbf{w}'_i \boldsymbol{\beta}_{10}) x_{ij},$$

which is the "covariance" between the $j$th covariate and the response variable. With the centering and standardization given in (2), $\xi_{nj}/\sigma$ is the correlation coefficient.

(B1) (a) $\varepsilon_i, \varepsilon_2, \ldots$ are independent and identically distributed random variables with mean zero and variance $\sigma^2$, where $0 < \sigma^2 < \infty$; (b) $\varepsilon_i$'s are sub-Gaussian, that is, their tail probabilities satisfy $P(|\varepsilon_i| > x) \leq K \exp(-Cx^2), i = 1, 2, \ldots$, for constants $C$ and $K$.

(B2) (Partial orthogonality) (a) There exists a constant $c_0 > 0$ such that

$$\left| n^{-1/2} \sum_{i=1}^{n} x_{ij} x_{ik} \right| \leq c_0, \qquad j \in J_n, k \in K_n,$$

for all $n$ sufficiently large. (b) There exists a constant $\xi_0 > 0$ such that $\min_{k \in K_n} |\xi_{nj}| > \xi_0 > 0$.

(B3) (a) $\lambda_n/n \to 0$ and $\lambda_n n^{-\gamma/2} k_n^{\gamma-2} \to \infty$; (b) $\log(m_n) = o(1) \times (\lambda_n n^{-\gamma/2})^{2/(2-\gamma)}$.

(B4) There exist constants $0 < b_1 < \infty$ such that $\max_{k \in K_n} |\beta_{1k}| \leq b_1$.

Condition (B1)(b) assumes that the tails of the error distribution behave like normal tails. Thus, it excludes heavy-tailed distributions. Condition (B2)(a) assumes that the covariates of the nonzero coefficients and the covariates of the zero coefficients are only weakly correlated. Condition (B2)(b) requires that the correlations between the covariates with nonzero coefficients and the dependent variable are bounded away from zero. Condition (B3)(a) restricts the penalty parameter $\lambda_n$ and the number of nonzero coefficients $k_n$. For $\lambda_n$, we must have $\lambda_n = o(n)$. For such a $\lambda_n$, $\lambda_n n^{-\gamma/2} k_n^{\gamma-2} =$ $o(1)n^{(2-\gamma)/2}k_n^{\gamma-2} = o(1)(n^{1/2}/k_n)^{2-\gamma}$. Thus, $k_n$ must satisfy $k_n/n^{1/2} = o(1)$. (B3)(b) restricts the number of zero coefficients $m_n$. To get a sense how large $m_n$ can be, we note that $\lambda_n$ can be as large as $\lambda_n = o(n)$. Thus, $\log(m_n) = o(1)(n^{(2-\gamma)/2})^{2/(2-\gamma)} = o(1)n$. So $m_n$ can be of the order $\exp(o(n))$. This certainly permits $m_n/n \to \infty$ and, hence, $p_n/n \to \infty$ as $n \to \infty$. Similar phenomena occur in Van der Laan and Bryan (2001) and Kosorok and Ma (2007) for uniformly consistent marginal estimators under different "large $p$, small $n$" data settings. On the other hand, the number of nonzero coefficients $k_n$ still must be smaller than $n$.

THEOREM 3. *Suppose that conditions* (B1) *to* (B4) *hold and that* $0 < \gamma < 1$. *Then*

$$P(\widetilde{\boldsymbol{\beta}}_{n2} = \mathbf{0}) \to 1 \quad \text{and} \quad P(\widetilde{\beta}_{n1k} \neq 0, k \in K_n) \to 1.$$

This theorem says that marginal bridge estimators can correctly distinguish between covariates with nonzero and zero coefficients with probability converging to one. However, the estimators of the nonzero coefficients are not consistent. To obtain consistent estimators, we use a two-step approach. First, we use the marginal bridge estimator to select the covariates with nonzero coefficients. Then we estimate the regression model with the selected covariates. In the second step, any reasonable regression method can be used. The choice of method is likely to depend on the characteristics of the data at hand, including the number of nonzero coefficients selected in the first step, the properties of the design matrix and the shape of the distribution of the $\varepsilon_i$'s. A two-step approach different from the one proposed here was also used by Bair et al. (2006) in their approach for supervised principal component analysis. In a recent paper Zhao and Yu (2006) provided an irrepresentable condition under which the LASSO is variable selection consistent. It would be interesting to study the implications of the irrepresentable condition in the context of bridge regression.

We now consider the use of the bridge objective function for second-stage estimation of $\boldsymbol{\beta}_{10}$, the vector of nonzero coefficients. Since the zero coefficients are correctly identified with probability converging to one, we can assume that only the covariates with nonzero coefficients are included in the model in the asymptotic analysis of the second step estimation. Let $\widehat{\boldsymbol{\beta}}_{1n}^*$ be the estimator in this step. Then, for the purpose of deriving its asymptotic distribution, it can be defined as the value that minimizes

$$(6) \quad U_n^*(\boldsymbol{\beta}_1) = \sum_{i=1}^{n}(y_i - \mathbf{w}_i'\boldsymbol{\beta}_1)^2 + \lambda_n^* \sum_{j=1}^{k_n}|\beta_{1j}|^\gamma,$$

where $\boldsymbol{\beta}_1 = (\beta_{11}, \ldots, \beta_{1k_n})'$. In addition to conditions (B1) to (B4), we assume the following:

(B5) (a) There exist a constant $\tau_1 > 0$ such that $\tau_{1n} \geq \tau_1$ for all $n$ sufficiently large;

(b) The covariates of nonzero coefficients satisfy $n^{-1/2}\max_{1\leq i \leq n} \mathbf{w}'_i\mathbf{w}_i \to 0$.

(B6) (a) $k_n(1+\lambda_n^*)/n \to 0$; (b) $\lambda_n^*(k_n/n)^{1/2} \to 0$.

These two conditions are needed for the asymptotic normality of $\widehat{\boldsymbol{\beta}}_{1n}^*$. Compared to condition (A5)(a), (B5)(a) assumes that the smallest eigenvalue of $\Sigma_{1n}$ is bounded away from zero, but does not assume that its largest eigenvalue is bounded. Condition (B5)(b) is the same as (A5)(b). In condition (B6), we can set $\lambda_n^* = 0$ for all $n$. Then $\widehat{\boldsymbol{\beta}}_{1n}^*$ is the OLS estimator. Thus, Theorem 4 below is applicable to the OLS estimator. When $\lambda_n^*$ is zero, then (B6)(a) becomes $k_n/n \to 0$ and (B6)(b) is satisfied for any value of $k_n$. Condition (B5)(b) also restricts $k_n$ implicitly. For example, if the covariates in $\mathbf{w}_i$ are bounded below by a constant $w_0 > 0$, then $\mathbf{w}'_i\mathbf{w}_i \geq k_n w_0^2$. So for (B5)(b) to hold, we must have $k_n n^{-1/2} \to 0$.

THEOREM 4. *Suppose that conditions* (B1) *to* (B6) *hold and that* $0 < \gamma < 1$. *Let* $s_n^2 = \sigma^2 \boldsymbol{\alpha}'_n \Sigma_{1n}^{-1}\boldsymbol{\alpha}_n$ *for any* $k_n \times 1$ *vector* $\boldsymbol{\alpha}_n$ *satisfying* $\|\boldsymbol{\alpha}_n\|_2 \leq 1$. *Then*

$$(7) \quad n^{1/2}s_n^{-1}\boldsymbol{\alpha}'_n(\widehat{\boldsymbol{\beta}}_{1n}^* - \boldsymbol{\beta}_{10}) = n^{-1/2}s_n^{-1}\sum_{i=1}^n \varepsilon_i \boldsymbol{\alpha}'_n \Sigma_{1n}^{-1}\mathbf{w}_i + o_p(1) \to_D N(0,1),$$

*where* $o_p(1)$ *is a term that converges to zero in probability uniformly with respect to* $\boldsymbol{\alpha}_n$.

**4. Numerical studies.** In this section we use simulation to evaluate the finite sample performance of bridge estimators.

4.1. *Computation of bridge estimators.* The penalized objective function (1) is not differentiable when $\beta$ has zero components. This singularity causes standard gradient based methods to fail. Motivated by the method of Fan and Li (2001) and Hunter and Li (2005), we approximate the bridge penalty by a function that has finite gradient at zero. Specifically, we approximate the bridge penalty function by $\sum_{j=1}^{p_n} \int_{-\infty}^{\beta_j} [\text{sgn}(u)/(|u|^{1/2} + \eta)]\,du$ for a small $\eta > 0$. We note this function and its gradient converge to the bridge penalty and its gradient as $\eta \to 0$, respectively.

Let $p = p_n$ be the number of covariates. Let $\widehat{\boldsymbol{\beta}}^{(m)}$ be the value of the $m$th iteration from the optimization algorithm, $m = 0, 1, \ldots$. Let $\tau$ be a prespecified convergence criterion. We set $\tau = 10^{-4}$ in our numerical studies. We conclude convergence if $\max_{1\leq j \leq p}|\widehat{\beta}_j^{(m)} - \widehat{\beta}_j^{(m+1)}| \leq \tau$, and conclude $\widehat{\beta}_j^{(m)} = 0$, if $|\widehat{\beta}_j^{(m)}| \leq \tau$. Denote $\mathbf{y}_n = (Y_1, \ldots, Y_n)$.

Initialize $\widehat{\boldsymbol{\beta}}^{(0)} = \mathbf{0}$ and $\eta = \tau$. For $m = 0, 1, \ldots$:

1. Compute the gradient of the sum of the squares $\mathbf{g}_1 = \mathbf{X}'_n(\mathbf{y}_n - \mathbf{X}_n\widehat{\boldsymbol{\beta}}^{(m)})$ and the approximate gradient of the penalty

   $$\mathbf{g}_2(\eta) = \tfrac{1}{2}\lambda(\operatorname{sgn}(\widehat{\beta}_1^{(m)})/(|\widehat{\beta}_1^{(m)}|^{1/2}+\eta),\ldots,\operatorname{sgn}(\widehat{\beta}_p^{(m)})/(|\widehat{\beta}_p^{(m)}|^{1/2}+\eta))'.$$

   Here $\mathbf{g}_1$ and $\mathbf{g}_2$ are $p\times 1$ vectors, with $j$th components $\mathbf{g}_{1j}$ and $\mathbf{g}_{2j}$, respectively. Note we use the notation $\mathbf{g}_2(\eta)$ to emphasize that the approximate gradient depends on $\eta$.
2. Compute the gradient $\mathbf{g}$ whose $j$th component, $g_j$, is defined as

   $$\text{if } |\widehat{\beta}_j^{(m)}| > \tau, \qquad g_j = g_{1j} + g_{2j}(\eta);$$
   $$\text{if } |\widehat{\beta}_j^{(m)}| \leq \tau, \qquad g_j = g_{1j} + g_{2j}(\eta^*),$$

   where $\eta^* = \arg\max_{j:0<|\widehat{\beta}_j^{(m)}|\leq\tau} |g_{1j}/g_{2j}(\eta)|$. In this way, we guarantee that, for the zero estimates, the corresponding components in $\mathbf{g}_2$ dominate the corresponding components in $\mathbf{g}_1$. Update $\eta = \eta^*$.
3. Re-scale $\mathbf{g} = \mathbf{g}/\max_j |g_j|$, such that its maximum component (in terms of absolute value) is less than or equal to 1. This step and the previous one guarantee that the increment in the components of $\boldsymbol{\beta}$ is less than $\tau$, the convergence criterion.
4. Update $\widehat{\boldsymbol{\beta}}^{(m+1)} = \widehat{\boldsymbol{\beta}}^{(m)} + \Delta \times \mathbf{g}$, where $\Delta$ is the increment in this iterative process. In our implementation we used $\Delta = 2 \times 10^{-3}$.
5. Replace $m$ by $m+1$ and repeat steps 1–5 until convergence.

Extensive simulation studies show that estimates obtained using this algorithm are well behaved and convergence is achieved under all simulated settings.

4.2. *Computation of marginal bridge estimators.* For a given penalty parameter $\lambda_n$, minimization of the marginal objective function $U_n$ defined in (4) amounts to solving a series of univariate minimization problems. Furthermore, since marginal bridge estimators are used only for variable selection, we do not need to solve the minimization problem. We only need to determine which coefficients are zero and which are not.

The objective function of each univariate minimization problem can be written in the form

$$g(u) = u^2 - 2au + \lambda|u|^\gamma,$$

where $|a| > 0$. By Lemma A of Knight and Fu (2000), $\arg\min(g) = 0$ if and only if

$$\lambda > \left(\frac{2}{2-\gamma}\right)\left(\frac{2(1-\gamma)}{2-\gamma}\right)^{1-\gamma}|a|^{2-\gamma}.$$

Therefore, computation for variable selection based on marginal bridge estimators can be done very quickly.

4.3. *Simulation studies.* This section describes simulation studies that are used to evaluate the finite sample performance of the bridge estimator. We investigate three features: (i) variable selection; (ii) prediction; and (iii) estimation. For (i), we measure variable selection performance by the frequency of correctly identifying zero and nonzero coefficients in repeated simulations. For (ii), we measure prediction performance using prediction mean square errors (PMSE), which are calculated from the fitted values based on the training data and the observed responses in an independent testing data not used in model fitting. For (iii), we measure estimation performance using the estimation mean square errors (EMSE) of the estimator, which are calculated from the estimated and true values of the parameters.

For comparison of prediction performance, we compare the PMSE of the bridge estimator to those of ordinary least squares (OLS) when applicable, ridge regression (RR), LASSO and Enet estimators. We assess the oracle property based on the variable selection results and the EMSE. For the bridge estimator, we set $\gamma = 1/2$. The RR, LASSO and elastic-net estimators are computed using the publicly available R packages (http://www.r-project.org). The bridge estimator is computed using the algorithm described in Section 4.1. The simulation scheme is close to the one in Zou and Hastie (2005), but differs in that the covariates are fixed instead of random.

We simulate data from the model

$$y = \mathbf{x}'\boldsymbol{\beta} + \epsilon, \qquad \epsilon \sim N(0, \sigma^2).$$

Six examples are considered, representing six different and commonly encountered scenarios. In each example the covariate vector $\mathbf{x}$ is generated from a multivariate normal distribution whose marginal distributions are standard $N(0,1)$ and whose covariance matrix is given in the description below. The value of $\mathbf{x}$ is generated once and then kept fixed. Replications are obtained by simulating the values of $\epsilon$ from $N(0, \sigma^2)$ and then setting $y = \mathbf{x}'\boldsymbol{\beta} + \epsilon$ for the fixed covariate value $\mathbf{x}$. Summary statistics are computed based on 500 replications. We consider six simulation models.

EXAMPLE 1.  $p = 30$ and $\sigma = 1.5$. The pairwise correlation between the $i$th and the $j$th components of $\mathbf{x}$ is $r^{|i-j|}$ with $r = 0.5$. Components 1–5 of $\boldsymbol{\beta}$ are 2.5; components 6–10 are 1.5; components 11–15 are 0.5 and the rest are zero. So there are 15 nonzero covariate effects five large effects, five moderate effects and five small effects.

EXAMPLE 2.  The same as Example 1, except that $r = 0.95$.

EXAMPLE 3. $p = 30$ and $\sigma = 1.5$. The predictors in Example 3 are generated as follows:

$$x_i = Z_1 + e_i, \qquad Z_1 \sim N(0,1), \qquad i = 1, \ldots, 5;$$
$$x_i = Z_2 + e_i, \qquad Z_2 \sim N(0,1), \qquad i = 6, \ldots, 10;$$
$$x_i = Z_3 + e_i, \qquad Z_3 \sim N(0,1), \qquad i = 11, \ldots, 15;$$
$$x_i \sim N(0,1), \qquad x_i \text{ i.i.d.} \qquad i = 16, \ldots, 30,$$

where $e_i$ are i.i.d. $N(0, 0.01), i = 1, \ldots, 15$. The first 15 components of $\beta$ are 1.5, the remaining ones are zero.

EXAMPLE 4. $p = 200$ and $\sigma = 1.5$. The first 15 covariates $(x_1, \ldots, x_{15})$ and the remaining 185 covariates $(x_{16}, \ldots, x_{200})$ are independent. The pairwise correlation between the $i$th and the $j$th components of $(x_1, \ldots, x_{15})$ is $r^{|i-j|}$ with $r = 0.5, i, j = 1, \ldots, 15$. The pairwise correlation between the $i$th and the $j$th components of $(x_{16}, \ldots, x_{200})$ is $r^{|i-j|}$ with $r = 0.5, i, j = 16, \ldots, 200$. Components 1–5 of $\beta$ are 2.5, components 6–10 are 1.5, components 11–15 are 0.5 and the rest are zero. So there are 15 nonzero covariate effects—five large effects, five moderate effects and five small effects. The covariate matrix has the partial orthogonal structure.

EXAMPLE 5. The same as Example 4, except that $r = 0.95$.

EXAMPLE 6. $p = 500$ and $\sigma = 1.5$. The first 15 covariates are generated the same way as in Example 5. The remaining 485 covariates are independent of the first 15 covariates and are generated independently from $N(0,1)$. The first 15 coefficients equal 1.5, and the remaining 485 coefficients are zero.

The examples with $r = 0.5$ have weak to moderate correlation among covariates, whereas those with $r = 0.95$ have moderate to strong correlations among covariates. Examples 3 and 6 correspond to the "grouping effects" in Zou and Hastie (2005) with three equally important groups. In Examples 3 and 6, covariates within the same group are highly correlated and the pairwise correlation coefficients are as high as 0.99. Therefore, there is particularly strong collinearity among the covariates in these two examples.

Following the simulation approach of Zou and Hastie (2005), in each example, the simulated data consist of a training set and an independent validation set and an independent test set, each of size 100. The tuning parameter is selected using the same simple approach as in Zou and Hastie (2005). We first fit the model with a given tuning parameter using the training set data only and compute the mean squared error between the fitted values and the responses in the validation data. We then search the tuning parameter space and choose the one with the smallest mean squared error as the final penalty parameter. Using this penalty parameter and the model estimated based on the training set, we compute the PMSE for the testing set. We also compute the probabilities that the estimators correctly identify covariates with nonzero and zero coefficients.

In Examples 1–3, the number of covariates is less than the sample size, so we use the bridge approach directly with the algorithm of Section 4.1. In Examples 4–6, the number of covariates is greater than the sample size. We use the two-step approach described in Section 3. We first select the nonzero covariates using the marginal bridge method. The number of nonzero covariates identified is much less than the sample size. In the second step, we use OLS.

Summary statistics of the variable selection and PMSE results based on 500 replicates are shown in Table 1. We see that the numbers of nonzero covariates selected by the bridge estimators are close to the true value (=15) in all examples. This agrees with the consistent variable selection result of Theorem 2. On average, the bridge estimator outperforms LASSO and ENet in terms of variable selection. Table 1 also gives the PMSEs of the Bridge, RR, LASSO and Enet estimators. For OLS (when applicable), LASSO, ENet and Bridge, the PMSEs are mainly caused by the variance of the random error. So the PMSEs are close, in general, with the Enet and Bridge being better than the LASSO and OLS. The RR is less satisfactory in Examples 4–6 with 200 covariates.

Figure 1 shows the frequencies of individual covariate effects being correctly "classified": zero versus nonzero. For better resolution, we only plot the first 30 covariates for Examples 4–6. We see that the bridge estimator can effectively identify large and moderate nonzero covariate effects and zero covariate effects.

TABLE 1
*Simulation study: comparison of OLS, RR, LASSO, Elastic net and the bridge estimator with $\gamma = 1/2$. PMSE: median of PMSE, inside "$(\cdot)$" are the corresponding standard deviations. Covariate: median of number of covariates with nonzero coefficients*

| Example |           | OLS         | RR          | LASSO       | ENet        | Bridge      |
|---------|-----------|-------------|-------------|-------------|-------------|-------------|
| 1       | PMSE      | 3.32 (0.58) | 3.51 (0.69) | 2.92 (0.51) | 2.80 (0.47) | 2.95 (0.51) |
|         | Covariate | 30          | 30          | 23          | 22          | 17          |
| 2       | PMSE      | 3.21 (0.53) | 2.65 (0.41) | 2.60 (0.40) | 2.46 (0.35) | 2.37 (0.36) |
|         | Covariate | 30          | 30          | 18          | 16          | 15          |
| 3       | PMSE      | 3.26 (0.58) | 3.34 (0.58) | 2.66 (0.40) | 2.38 (0.33) | 2.31 (0.34) |
|         | Covariate | 30          | 30          | 18          | 15          | 15          |
| 4       | PMSE      | –           | 20.45 (2.02)| 3.55 (0.64) | 3.30 (0.53) | 3.98 (0.83) |
|         | Covariate | –           | 200         | 37          | 37          | 29          |
| 5       | PMSE      | –           | 5.80 (1.31) | 2.71 (0.42) | 2.50 (0.36) | 2.64 (0.44) |
|         | Covariate | –           | 200         | 25          | 16          | 15          |
| 6       | PMSE      | –           | 43.10 (2.23)| 3.51 (0.57) | 2.70 (0.49) | 2.68 (0.39) |
|         | Covariate | –           | 500         | 43          | 20          | 17          |

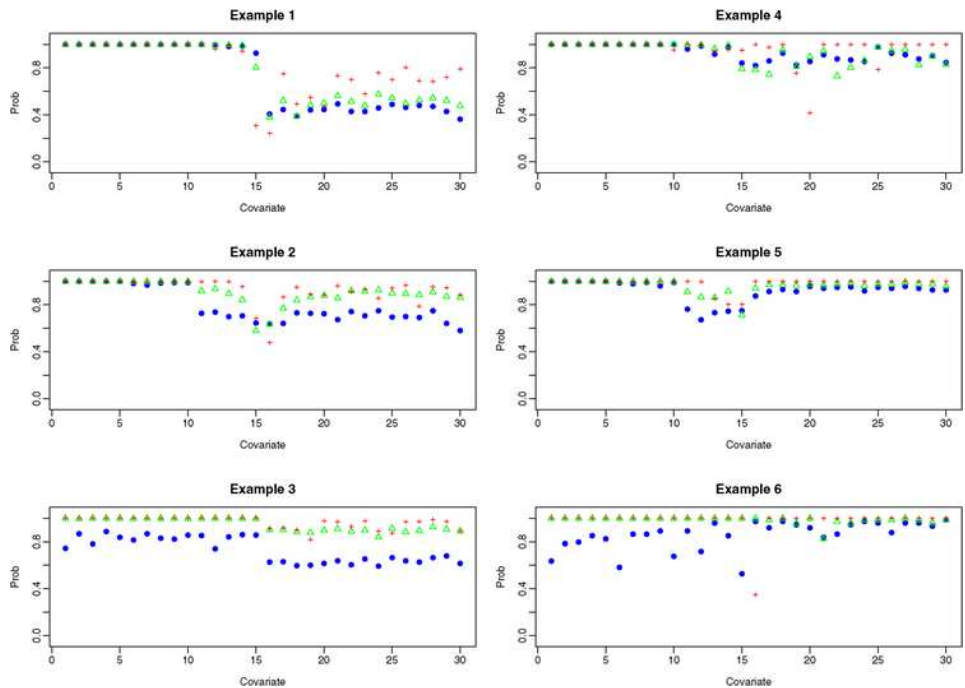

FIG. 1. *Simulation study (Examples 1–6): probability of individual covariate effect being correctly identified. Circle: LASSO; Triangle: ENet; Plus sign: Bridge estimate.*

TABLE 2
*Simulation study: comparison of OLS with the first 15 covariates (OLS-oracle), bridge estimate with the first 15 covariates (bridge-oracle) and bridge estimate with all covariates. For each model, the first row: median of absolute bias (across the 15 covariates) and median of variance (across the 15 covariates); the second row: median of EMSE and standard deviation of EMSE*

| Example | | OLS-oracle | Bridge-oracle | Bridge |
|---|---|---|---|---|
| 1 | bias/sd | 0.007, 0.047 | 0.019, 0.045 | 0.035, 0.020 |
|   | EMSE    | 0.647, 0.306 | 0.625, 0.305 | 0.702, 0.311 |
| 2 | bias/sd | 0.014, 0.509 | 0.114, 0.053 | 0.024, 0.018 |
|   | EMSE    | 7.252, 3.707 | 0.910, 1.109 | 0.990, 0.738 |
| 3 | bias/sd | 0.041, 2.041 | 0.026, 0.080 | 0.028, 0.007 |
|   | EMSE    | 30.15, 14.01 | 0.163, 3.468 | 0.133, 0.898 |
| 4 | bias/sd | 0.006, 0.043 | 0.014, 0.042 | 0.061, 0.062 |
|   | EMSE    | 0.655, 0.293 | 0.662, 0.281 | 1.186, 0.849 |
| 5 | bias/sd | 0.036, 0.535 | 0.133, 0.051 | 0.050, 0.467 |
|   | EMSE    | 7.077, 3.565 | 1.179, 0.714 | 7.013, 3.629 |
| 6 | bias/sd | 0.035, 1.928 | 0.027, 0.078 | 0.072, 1.923 |
|   | EMSE    | 28.90, 12.46 | 0.218, 2.967 | 28.43, 12.65 |

Simulation studies were also carried out to investigate the asymptotic oracle property of the bridge estimator. This property says that bridge estimators have the same asymptotic efficiency as the estimator obtained under the knowledge of which coefficients are nonzero and which are zero. To evaluate this property, we consider three estimators: OLS using the covariates with nonzero coefficients only (OLS-oracle); the bridge estimator using the covariates with nonzero coefficients (bridge-oracle); and the bridge estimator using all the covariates. We note that the OLS-oracle and bridge-oracle estimators cannot be used in practice. We use them here only for the purpose of comparison. We use the same six examples as described above.

Table 2 presents the summary statistics based on 500 replications. In Examples 1–3, the bridge estimator and bridge-oracle estimators perform similarly. In Examples 4–6, the bridge estimator is similar to the OLS-oracle estimator. In Examples 2 and 3 where the covariates are highly correlated, the OLS-oracle estimators have considerably larger EMSEs than the bridge-oracle and bridge estimators. In Examples 4 and 6, the OLS-oracle estimators and the two-step estimators have considerably larger EMSEs than the bridge-oracle estimators. This is due to the fact that OLS estimators tend to perform poorly when there is strong collinearity among covariates. The simulation results from these examples also suggest that, in finite samples, bridge estimators provide substantial improvement over the OLS estimators in terms of EMSE in the presence of strong collinearity.

**5. Concluding remarks.** In this paper we have studied the asymptotic properties of bridge estimators when the number of covariates and regression coefficients increases to infinity as $n \to \infty$. We have shown that, when $0 < \gamma < 1$, bridge estimators correctly identify zero coefficients with probability converging to one and that the estimators of nonzero coefficients are asymptotically normal and oracle efficient. Our results generalize the results of Knight and Fu (2000), who studied the asymptotic behavior of LASSO-type estimators in the finite-dimensional regression parameter setting. Theorems 1 and 2 were obtained under the assumption that the number of parameters is smaller than the sample size, as described in conditions (A2) and (A3). They are not applicable when the number of parameters is greater than the sample size, which arises in microarray gene expression studies. Accordingly, we have also considered a marginal bridge estimator under a partial orthogonality condition in which the covariates of zero coefficients are orthogonal to or only weakly correlated with the covariates of nonzero coefficients. The marginal bridge estimator can consistently distinguish covariates with zero and nonzero coefficients even when the number of zero coefficients is greater than the sample size. Indeed, the number of zero coefficients can be in the order of $\exp(o(n))$.

We have proposed a gradient based algorithm for computing bridge estimators. Our simulation study suggests this algorithm converges reasonably rapidly. It also suggests that the bridge estimator with $\gamma = 1/2$ behaves well in our simulated models. The bridge estimator correctly identifies zero coefficients with higher probability than do the LASSO and Elastic-net estimators. It also performs well in terms of predictive mean square errors. Our theoretical and numerical results suggest that the bridge estimator with $0 < \gamma < 1$ is a useful alternative to the existing methods for variable selection and parameter estimation with high-dimensional data.

**6. Proofs.** In this section we give the proofs of the results stated in Sections 2 and 3. For simplicity of notation and without causing confusion, we write $\mathbf{X}_n$, $\mathbf{X}_{1n}$ and $\mathbf{X}_{2n}$ as $\mathbf{X}$, $\mathbf{X}_1$ and $\mathbf{X}_2$.

We first prove the following lemma which will be used in the proof of Theorem 1.

LEMMA 1. *Let $\mathbf{u}$ be a $p_n \times 1$ vector. Under condition* (A1)(a),
$$\operatorname{E} \sup_{\|u\| < \delta} \left| \sum_{i=1}^n \varepsilon_i \mathbf{x}_i' \mathbf{u} \right| \leq \delta \sigma n^{1/2} p_n^{1/2}.$$

PROOF. By the Cauchy–Schwarz inequality and condition (A1), we have
$$\operatorname{E} \sup_{\|\mathbf{u}\| \leq \delta} \left| \sum_{i=1}^n \varepsilon_i \mathbf{x}_i' \mathbf{u} \right|^2 \leq \operatorname{E} \sup_{\|\mathbf{u}\| \leq \delta} \|\mathbf{u}\|^2 \left\| \sum_{i=1}^n \varepsilon_i \mathbf{x}_i \right\|^2$$
$$\leq \delta^2 \operatorname{E} \left[ \sum_{i=1}^n \varepsilon_i \mathbf{x}_i' \sum_{i=1}^n \varepsilon_i \mathbf{x}_i \right]$$
$$= \delta^2 \sigma^2 \sum_{i=1}^n \mathbf{x}_i' \mathbf{x}_i$$
$$= \delta^2 \sigma^2 n \operatorname{trace}\left( n^{-1} \sum_{i=1}^n \mathbf{x}_i \mathbf{x}_i' \right)$$
$$= \delta^2 \sigma^2 n p_n.$$

Thus, the lemma follows from Jensen's inequality. □

PROOF OF THEOREM 1. We first show that
$$\|\widehat{\boldsymbol{\beta}}_n - \boldsymbol{\beta}_0\| = O_p((p_n + \lambda_n k_n)/(n\rho_{1n}))^{1/2}. \tag{8}$$

By the definition of $\widehat{\boldsymbol{\beta}}_n$,

$$\sum_{i=1}^{n}(Y_i - \mathbf{x}_i'\widehat{\boldsymbol{\beta}}_n)^2 + \lambda_n \sum_{j=1}^{p_n} |\widehat{\beta}_j|^\gamma \leq \sum_{i=1}^{n}(Y_i - \mathbf{x}_i'\boldsymbol{\beta}_0)^2 + \lambda_n \sum_{j=1}^{p_n} |\beta_{0j}|^\gamma.$$

It follows that

$$\sum_{i=1}^{n}(Y_i - \mathbf{x}_i'\widehat{\boldsymbol{\beta}}_n)^2 \leq \sum_{i=1}^{n}(Y_i - \mathbf{x}_i'\boldsymbol{\beta}_0)^2 + \lambda_n \sum_{j=1}^{p_n} |\beta_{0j}|^\gamma.$$

Let $\eta_n = \lambda_n \sum_{j=1}^{p_n} |\beta_{0j}|^\gamma$, then

$$\eta_n \geq \sum_{i=1}^{n}(Y_i - \mathbf{x}_i'\widehat{\boldsymbol{\beta}}_n)^2 - \sum_{i=1}^{n}(Y_i - \mathbf{x}_i'\boldsymbol{\beta}_0)^2$$

$$= \sum_{i=1}^{n}[\mathbf{x}_i'(\widehat{\boldsymbol{\beta}}_n - \boldsymbol{\beta}_0)]^2 + 2\sum_{i=1}^{n} \varepsilon_i \mathbf{x}_i(\boldsymbol{\beta}_0 - \widehat{\boldsymbol{\beta}}_n).$$

Let $\boldsymbol{\delta}_n = n^{1/2}(\Sigma_n)^{1/2}(\widehat{\boldsymbol{\beta}}_n - \boldsymbol{\beta}_0)$, $\mathbf{D}_n = n^{-1/2}(\Sigma_n)^{-1/2}\mathbf{X}'$ and $\boldsymbol{\varepsilon}_n = (\varepsilon_1, \ldots, \varepsilon_n)'$. Then

$$\sum_{i=1}^{n}[\mathbf{x}_i'(\widehat{\boldsymbol{\beta}}_n - \boldsymbol{\beta}_0)]^2 + 2\sum_{i=1}^{n} \varepsilon_i \mathbf{x}_i'(\boldsymbol{\beta}_0 - \widehat{\boldsymbol{\beta}}_n) = \boldsymbol{\delta}_n'\boldsymbol{\delta}_n - 2(\mathbf{D}_n\boldsymbol{\varepsilon}_n)'\boldsymbol{\delta}_n.$$

So we have $\boldsymbol{\delta}_n'\boldsymbol{\delta}_n - 2(\mathbf{D}_n\boldsymbol{\varepsilon})'\boldsymbol{\delta}_n - \eta_n \leq 0$. That is, $\|\boldsymbol{\delta}_n - \mathbf{D}_n\boldsymbol{\varepsilon}_n\|^2 - \|\mathbf{D}_n\boldsymbol{\varepsilon}_n\|^2 - \eta_n \leq 0$. Therefore, $\|\boldsymbol{\delta}_n - \mathbf{D}_n\boldsymbol{\varepsilon}_n\| \leq \|\mathbf{D}_n\boldsymbol{\varepsilon}_n\| + \eta_n^{1/2}$. By the triangle inequality,

$$\|\boldsymbol{\delta}_n\| \leq \|\boldsymbol{\delta}_n - \mathbf{D}_n\boldsymbol{\varepsilon}_n\| + \|\mathbf{D}_n\boldsymbol{\varepsilon}_n\| \leq 2\|\mathbf{D}_n\boldsymbol{\varepsilon}_n\| + \eta_n^{1/2}.$$

It follows that $\|\boldsymbol{\delta}_n\|^2 \leq 6\|\mathbf{D}_n\boldsymbol{\varepsilon}_n\|^2 + 3\eta_n$. Let $\mathbf{d}_i$ be the $i$th column of $\mathbf{D}_n$. Then $\mathbf{D}_n\boldsymbol{\varepsilon} = \sum_{i=1}^n \mathbf{d}_i\varepsilon_i$. Since $\mathrm{E}\varepsilon_i\varepsilon_j = 0$ if $i \neq j$, $\mathrm{E}\|\mathbf{D}_n\boldsymbol{\varepsilon}_n\|^2 = \sum_{i=1}^n \|\mathbf{d}_i\|^2 \mathrm{E}\varepsilon_i^2 = \sigma^2 \mathrm{tr}(\mathbf{D}_n\mathbf{D}_n') = \sigma^2 p_n$. So we have $\mathrm{E}\|\boldsymbol{\delta}_n\|^2 \leq 6\sigma^2 p_n + 3\eta_n$. That is,

(9) $$n\mathrm{E}[(\widehat{\boldsymbol{\beta}}_n - \boldsymbol{\beta}_0)'\Sigma_n(\widehat{\boldsymbol{\beta}}_n - \boldsymbol{\beta}_0)] \leq 6\sigma^2 p_n + 3\eta_n.$$

Since the number of nonzero coefficients is $k_n$, $\eta_n = \lambda_n \sum_{j=1}^{p_n} |\beta_{0j}|^\gamma = O(\lambda_n k_n)$. Noting that $\rho_{1n}$ is the smallest eigenvalue of $\Sigma_{1n}$, (8) follows from (9).

We now show that

(10) $$\|\widehat{\boldsymbol{\beta}}_n - \boldsymbol{\beta}_0\| = O_p(\rho_{1n}^{-1}(p_n/n)^{1/2}).$$

Let $r_n = \rho_{1n}(n/p_n)^{1/2}$. The proof of (10) follows that of Theorem 3.2.5 of Van der Vaart and Wellner (1996). For each $n$, partition the parameter space (minus $\boldsymbol{\beta}_0$) into the "shells" $S_{j,n} = \{\boldsymbol{\beta} : 2^{j-1} < r_n\|\boldsymbol{\beta} - \boldsymbol{\beta}_0\| < 2^j\}$ with $j$ ranging over the integers. If $r_n\|\widehat{\boldsymbol{\beta}}_n - \boldsymbol{\beta}_0\|$ is larger than $2^M$ for a given

integer $M$, then $\widehat{\boldsymbol{\beta}}_n$ is in one of the shells with $j \geq M$. By the definition of $\widehat{\boldsymbol{\beta}}_n$ that it minimizes $L_n(\boldsymbol{\beta})$, for every $\epsilon > 0$,

$$\mathrm{P}(r_n\|\widehat{\boldsymbol{\beta}}_n - \boldsymbol{\beta}_0\| > 2^M)$$
$$= \sum_{j \geq M, 2^j \leq \epsilon r_n} \mathrm{P}\left(\inf_{\boldsymbol{\beta} \in S_{j,n}} (L_n(\boldsymbol{\beta}) - L_n(\boldsymbol{\beta}_0)) \leq 0\right) + \mathrm{P}(2\|\widehat{\boldsymbol{\beta}}_n - \boldsymbol{\beta}_0\| \geq \epsilon).$$

Because $\widehat{\boldsymbol{\beta}}_n$ is consistent by (8) and condition (A2), the second term on the right-hand side converges to zero. So we only need to show that the first term on the right-hand side converges to zero. Now

$$L_n(\boldsymbol{\beta}) - L_n(\boldsymbol{\beta}_0)$$
$$= \sum_{i=1}^{n}(Y_i - \mathbf{x}_i'\boldsymbol{\beta})^2 + \lambda_n \sum_{j=1}^{k_n}|\beta_{1j}|^\gamma + \lambda_n \sum_{j=1}^{m_n}|\beta_{2j}|^\gamma$$
$$- \sum_{i=1}^{n}(Y_i - \mathbf{w}_i'\boldsymbol{\beta}_{10})^2 - \lambda_n \sum_{i=1}^{k_n}|\beta_{01j}|^\gamma$$
$$\geq \sum_{i=1}^{n}(Y_i - \mathbf{x}_i'\boldsymbol{\beta})^2 + \lambda_n \sum_{j=1}^{k_n}|\beta_{1j}|^\gamma - \sum_{i=1}^{n}(Y_i - \mathbf{w}_i'\boldsymbol{\beta}_{10})^2 - \lambda_n \sum_{i=1}^{k_n}|\beta_{01j}|^\gamma$$
$$= \sum_{i=1}^{n}[\mathbf{x}_i'(\boldsymbol{\beta} - \boldsymbol{\beta}_0)]^2 - 2\sum_{i=1}^{n}\varepsilon_i \mathbf{x}_i'(\boldsymbol{\beta} - \boldsymbol{\beta}_0) + \lambda_n \sum_{j=1}^{k_n}\{|\beta_{1j}|^\gamma - |\beta_{01j}|^\gamma\}$$
$$\equiv I_{1n} + I_{2n} + I_{3n}.$$

On $S_{j,n}$, the first term $I_{1n} \geq n\rho_{1n}2^{2(j-1)}r_n^{-2}$. The third term

$$I_{3n} = \lambda_n \gamma \sum_{j=1}^{k_n} |\beta_{01j}^*|^{\gamma-1} \operatorname{sgn}(\beta_{01j})(\beta_{1j} - \beta_{01j}),$$

for some $\beta_{01j}^*$ between $\beta_{01j}$ and $\beta_{1j}$. By condition (A4) and since we only need to consider $\boldsymbol{\beta}$ with $\|\boldsymbol{\beta} - \boldsymbol{\beta}_0\| \leq \epsilon$, there exists a constant $c_3 > 0$ such that

$$|I_{3n}| \leq c_3 \gamma \lambda_n \sum_{j=1}^{k_n} |\beta_{1j} - \beta_{01j}| \leq c_3 \gamma \lambda_n k_n^{1/2} \|\boldsymbol{\beta} - \boldsymbol{\beta}_0\|.$$

So on $S_{j,n}$, $I_{3n} \geq -c_3 \lambda_n k_n^{1/2}(2^j/r_n)$. Therefore, on $S_{j,n}$,

$$L_n(\boldsymbol{\beta}) - L_n(\boldsymbol{\beta}_0) \geq -|I_{2n}| + n\rho_{1n}(2^{2(j-1)}/r_n^2) - c_3 \lambda_n k_n^{1/2}(2^j/r_n).$$

It follows that

$$P\left(\inf_{\boldsymbol{\beta}\in S_{j,n}}(L_n(\boldsymbol{\beta})-L_n(\boldsymbol{\beta}_0))\leq 0\right)$$
$$\leq P\left(\sup_{\boldsymbol{\beta}\in S_{j,n}}|I_{2n}|\geq n\rho_{1n}(2^{2(j-1)}/r_n^2)-c_3\lambda_n k_n^{1/2}(2^j/r_n)\right)$$
$$\leq \frac{2n^{1/2}p_n^{1/2}(2^j/r_n)}{n\rho_{1n}(2^{2(j-1)}/r_n^2)-c_3\lambda_n k_n^{1/2}(2^j/r_n)}$$
$$=\frac{2}{2^{j-2}-c_3\lambda_n k_n^{1/2}(np_n)^{-1/2}},$$

where the second inequality follows from Markov's inequality and Lemma 1. Under condition (A3)(a), $\lambda_n k_n^{1/2}(np_n)^{-1/2}\to 0$ as $n\to\infty$. So for $n$ sufficiently large, $2^{j-2}-c_3\lambda_n k_n^{1/2}(np_n)^{-1/2}\geq 2^{j-3}$ for all $j\geq 3$. Therefore,

$$\sum_{j\geq M, 2^j\leq \epsilon r_n}P\left(\inf_{\boldsymbol{\beta}\in S_{j,n}}(L_n(\boldsymbol{\beta})-L_n(\boldsymbol{\beta}_0))\leq 0\right)\leq \sum_{j\geq M}\frac{1}{2^{j-2}}\leq 2^{-(M-3)},$$

which converges to zero for every $M=M_n\to\infty$. This completes the proof of (10). Combining (8) and (10), the result follows. This completes the proof of Theorem 1. □

LEMMA 2. *Suppose that $0<\gamma<1$. Let $\widehat{\boldsymbol{\beta}}_n=(\widehat{\boldsymbol{\beta}}'_{1n},\widehat{\boldsymbol{\beta}}'_{2n})'$. Under conditions* (A1) *to* (A4), *$\widehat{\boldsymbol{\beta}}_{2n}=0$ with probability converging to 1.*

PROOF. By Theorem 1, for a sufficiently large $C$, $\widehat{\boldsymbol{\beta}}_n$ lies in the ball $\{\boldsymbol{\beta}:\|\boldsymbol{\beta}-\boldsymbol{\beta}_0\|\leq h_nC\}$ with probability converging to 1, where $h_n=\rho_{1n}^{-1}(p_n/n)^{1/2}$. Let $\boldsymbol{\beta}_{1n}=\boldsymbol{\beta}_{01}+h_n\mathbf{u}_1$ and $\boldsymbol{\beta}_{2n}=\boldsymbol{\beta}_{02}+h_n\mathbf{u}_2=h_n\mathbf{u}_2$ with $\|\mathbf{u}\|_2^2=\|\mathbf{u}_1\|_2^2+\|\mathbf{u}_2\|_2^2\leq C^2$. Let

$$V_n(\mathbf{u}_1,\mathbf{u}_2)=L_n(\boldsymbol{\beta}_{1n},\boldsymbol{\beta}_{2n})-L_n(\boldsymbol{\beta}_{10},\mathbf{0})=L_n(\boldsymbol{\beta}_{10}+h_n\mathbf{u}_1,h_n\mathbf{u}_2)-L_n(\boldsymbol{\beta}_{10},\mathbf{0}).$$

Then $\widehat{\boldsymbol{\beta}}_{1n}$ and $\widehat{\boldsymbol{\beta}}_{2n}$ can be obtained by minimizing $V_n(\mathbf{u}_1,\mathbf{u}_2)$ over $\|\mathbf{u}\|\leq C$, except on an event with probability converging to zero. To prove the lemma, it suffices to show that, for any $\mathbf{u}_1$ and $\mathbf{u}_2$ with $\|\mathbf{u}\|\leq C$, if $\|\mathbf{u}_2\|>0$, $V_n(\mathbf{u}_1,\mathbf{u}_2)-V_n(\mathbf{u}_1,\mathbf{0})>0$ with probability converging to 1. Some simple calculation shows that

$$V_n(\mathbf{u}_1,\mathbf{u}_2)-V_n(\mathbf{u}_1,0)=h_n^2\sum_{i=1}^n(\mathbf{z}'_i\mathbf{u}_2)^2+2h_n^2\sum_{i=1}^n(\mathbf{w}'_i\mathbf{u}_1)(\mathbf{z}'_i\mathbf{u}_2)$$
$$-2h_n\sum_{i=1}^n\varepsilon_i(\mathbf{z}'_i\mathbf{u}_2)+\lambda_n h_n^\gamma\sum_{j=1}^{m_n}|u_{2j}|^\gamma$$
$$\equiv II_{1n}+II_{2n}+II_{3n}+II_{4n}.$$

For the first two terms, we have

$$II_{1n} + II_{2n} \geq h_n^2 \sum_{i=1}^n (\mathbf{z}_i'\mathbf{u}_2)^2 - h_n^2 \sum_{i=1}^n [(\mathbf{w}_i'\mathbf{u}_1)^2 + (\mathbf{z}_i'\mathbf{u}_2)^2]$$

$$= -h_n^2 \sum_{i=1}^n (\mathbf{w}_i'\mathbf{u}_1)^2$$

(11)
$$\geq -nh_n^2 \tau_{2n} \|\mathbf{u}_1\|^2$$
$$\geq -\tau_2 (p_n/\rho_{1n}^2) C^2,$$

where we used condition (A5)(a) in the last inequality. For the third term, since

$$\mathrm{E}\left|\sum_{i=1}^n \varepsilon_i \mathbf{z}_i'\mathbf{u}_2\right| \leq \left[\mathrm{E}\left(\sum_{i=1}^n \varepsilon_i \mathbf{z}_i'\mathbf{u}_2\right)^2\right]^{1/2}$$

$$= \sigma \left[\sum_{i=1}^n \mathbf{u}_2' \mathbf{z}_i \mathbf{z}_i' \mathbf{u}_2\right]^{1/2}$$

$$\leq \sigma n^{1/2} \rho_{2n}^{1/2} \|\mathbf{u}_2\|$$

$$\leq \sigma (np_n)^{1/2} C,$$

we have

(12) $$II_{3n} = h_n n^{1/2} p_n^{1/2} O_p(1) = (p_n/\rho_{1n}) O_p(1).$$

For the fourth term, we first note that

$$\left[\sum_{j=1}^{m_n} |u_{2j}|^\gamma\right]^{2/\gamma} \geq \sum_{j=1}^{m_n} |u_{2j}|^2 = \|\mathbf{u}_2\|^2.$$

Thus,

(13) $$II_{4n} = \lambda_n h_n^\gamma O(\|\mathbf{u}_2\|^\gamma).$$

Under condition (A3)(b), $\lambda_n h_n^\gamma/(p_n \rho_{1n}^{-2}) = \lambda_n n^{\gamma/2} (\rho_{1n}/\sqrt{p_n})^{2-\gamma} \to \infty$. Combining (11), (12) and (13), we have, for $\|\mathbf{u}_2\|_2 > 0$, $V_n(\mathbf{u}) > 0$ with probability converging to 1. This completes the proof of Lemma 2. $\square$

PROOF OF THEOREM 2. Part (i) follows from Lemma 1. We need to prove (ii). Under conditions (A1) and (A2), $\widehat{\boldsymbol{\beta}}_n$ is consistent by Theorem 1. By condition (A4), each component of $\widehat{\boldsymbol{\beta}}_{1n}$ stays away from zero for $n$ sufficiently large. Thus, it satisfies the stationary equation evaluated at $(\widehat{\boldsymbol{\beta}}_{1n}, \widehat{\boldsymbol{\beta}}_{2n})$, $(\partial/\partial \boldsymbol{\beta}_1) L_n(\widehat{\boldsymbol{\beta}}_{1n}, \widehat{\boldsymbol{\beta}}_{2n}) = 0$. That is, $-2\sum_{i=1}^n (Y_i - \mathbf{w}_i' \widehat{\boldsymbol{\beta}}_{1n} -$ $\mathbf{z}_i'\widehat{\boldsymbol{\beta}}_{2n})\mathbf{w}_i + \lambda_n \gamma \psi_n = 0$, where $\psi_n$ is a $k_n \times 1$ vector whose $j$th element is $|\widehat{\beta}_{1nj}|^{\gamma-1}\operatorname{sgn}(\widehat{\beta}_{1nj})$. Since $\boldsymbol{\beta}_{20} = \mathbf{0}$ and $\varepsilon_i = Y_i - \mathbf{w}_i'\boldsymbol{\beta}_{10}$, this equation can be written $-2\sum_{i=1}^{n}(\varepsilon_i - \mathbf{w}_i'(\widehat{\boldsymbol{\beta}}_{1n} - \boldsymbol{\beta}_{10}) - \mathbf{z}_i'\widehat{\boldsymbol{\beta}}_{2n})\mathbf{w}_i + \lambda_n \gamma \psi_n = 0$. Therefore,

$$\Sigma_{1n}(\widehat{\boldsymbol{\beta}}_{1n} - \boldsymbol{\beta}_{10}) = \frac{1}{n}\sum_{i=1}^{n}\varepsilon_i \mathbf{w}_i - \frac{1}{2n}\gamma\lambda_n\psi_n - \frac{1}{n}\sum_{i=1}^{n}\mathbf{z}_i'\widehat{\boldsymbol{\beta}}_{2n}\mathbf{w}_i.$$

It follows that

$$n^{1/2}\boldsymbol{\alpha}_n'(\widehat{\boldsymbol{\beta}}_{1n} - \boldsymbol{\beta}_{10})$$
$$= n^{-1/2}\sum_{i=1}^{n}\varepsilon_i \boldsymbol{\alpha}_n'\Sigma_{1n}^{-1}\mathbf{w}_i - \tfrac{1}{2}\gamma n^{-1/2}\lambda_n\boldsymbol{\alpha}_n'\Sigma_{1n}^{-1}\psi_n - n^{-1/2}\sum_{i=1}^{n}\mathbf{z}_i'\widehat{\boldsymbol{\beta}}_{2n}\mathbf{w}_i.$$

By (i), $P(\widehat{\boldsymbol{\beta}}_{2n} = 0) \to 1$. Thus, the last term on the right-hand side equals zero with probability converging to 1. It certainly follows that it converges to zero in probability. When $\|\boldsymbol{\alpha}_n\| \leq 1$, under condition (A4),

$$|n^{-1/2}\boldsymbol{\alpha}_n'\Sigma_{1n}^{-1}\psi_n| \leq n^{-1/2}\zeta_1^{-1}\|\boldsymbol{\alpha}_n\| \cdot \|\|\widehat{\boldsymbol{\beta}}_{1n}|^{-(1-\gamma)}\|$$
$$\leq 2n^{-1/2}\tau_1^{-1}k_n^{1/2}b_0^{-(1-\gamma)},$$

except on an event with probability converging to zero. Under (A3)(a), $\lambda_n(k_n/n)^{1/2} \to 0$. Therefore,

(14) $\quad n^{1/2}s_n^{-1}\boldsymbol{\alpha}_n'(\widehat{\boldsymbol{\beta}}_{1n} - \boldsymbol{\beta}_{10}) = n^{-1/2}s_n^{-1}\sum_{i=1}^{n}\varepsilon_i\boldsymbol{\alpha}_n'\Sigma_{1n}^{-1}\mathbf{w}_i + o_p(1).$

We verify the conditions of the Lindeberg–Feller central limit theorem.

Let $v_i = n^{-1/2}s_n^{-1}\boldsymbol{\alpha}_n'\Sigma_{1n}^{-1}\mathbf{w}_i$ and $w_i = \varepsilon_i v_i$. First,

$$\operatorname{Var}\left(\sum_{i=1}^{n}w_i\right) = n^{-1}\sigma^2 s_n^{-2}\sum_{i=1}^{n}\boldsymbol{\alpha}_n'\Sigma_{1n}^{-1}\mathbf{w}_i\mathbf{w}_i'\Sigma_{1n}^{-1}\boldsymbol{\alpha}_n = s_n^{-2}s_n^2 = 1.$$

For any $\epsilon > 0$, $\sum_{i=1}^{n}E[w_i^2 \mathbf{1}\{|w_i| > \epsilon\}] = \sigma^2\sum_{i=1}^{n}v_i^2 E\varepsilon_i^2 \mathbf{1}\{|\varepsilon_i v_i| > \epsilon\}$. Since

$$\sigma^2\sum_{i=1}^{n}v_i^2 = n^{-1}\sigma^2 s_n^{-2}\sum_{i=1}^{n}(\boldsymbol{\alpha}_n'\Sigma_{1n}^{-1}\mathbf{w}_i\mathbf{w}_i'\Sigma_{1n}^{-1}\boldsymbol{\alpha}_n) = 1,$$

it suffices to show that, $\max_{1\leq i\leq n} E\varepsilon_i^2 \mathbf{1}\{|\varepsilon_i v_i| > \epsilon\} \to 0$, or equivalently,

(15) $\quad \max_{1\leq i\leq n}|v_i| = n^{-1/2}s_n^{-1}\max_{1\leq i\leq n}|\boldsymbol{\alpha}_n'\Sigma_{1n}^{-1}\mathbf{w}_i| \to 0.$

Since $|\boldsymbol{\alpha}_n'\Sigma_{1n}^{-1}\mathbf{w}_i| \leq (\boldsymbol{\alpha}_n'\Sigma_{1n}^{-1}\boldsymbol{\alpha}_n)^{1/2}(\mathbf{w}_i'\Sigma_{1n}^{-1}\mathbf{w}_i)^{1/2}$ and $s_n^{-1} = \sigma^{-1}(\boldsymbol{\alpha}_n\Sigma_{1n}^{-1}\boldsymbol{\alpha}_n)^{-1/2}$, we have

$$\max_{1\leq i\leq n}|v_i| \leq \sigma^{-1}n^{-1/2}\max_{1\leq i\leq n}(\mathbf{w}_i'\Sigma_{1n}^{-1}\mathbf{w}_i)^{1/2} \leq \sigma^{-1}\tau_1^{-1/2}n^{-1/2}\max_{1\leq i\leq n}(\mathbf{w}_i'\mathbf{w}_i)^{1/2},$$

(15) follows from assumption (A5). This completes the proof of Theorem 2.
□

LEMMA 3 [Knight and Fu (2000)]. *Let $g(u) = u^2 - 2au + \lambda|u|^\gamma$, where $a \neq 0$, $\lambda \geq 0$, and $0 < \gamma < 1$. Denote*

$$c_\gamma = \left(\frac{2}{2-\gamma}\right)\left(\frac{2(1-\gamma)}{2-\gamma}\right)^{1-\gamma}.$$

*Suppose that $a \neq 0$. Then $\arg\min(g) = 0$ if and only if $\lambda > c_\gamma |a|^{2-\gamma}$.*

Let $\psi_2(x) = \exp(x^2) - 1$. For any random variable $X$, its $\psi_2$-Orlicz norm $\|X\|_{\psi_2}$ is defined as $\|X\|_{\psi_2} = \inf\{C > 0 : \mathrm{E}\psi_2(|X|/C) \leq 1\}$. The Orlicz norm is useful for obtaining maximal inequalities; see Van der Vaart and Wellner (1996), Section 2.2.

LEMMA 4. *Let $c_1, \ldots, c_n$ be constants satisfying $\sum_{i=1}^n c_i^2 = 1$, and let $W = \sum_{i=1}^n c_i \varepsilon_i$.*

(i) *Under condition* (B1), *$\|W\|_{\psi_2} \leq K_2[\sigma + ((1+K)C^{-1})^{1/2}]$, where $K_2$ is a constant.*

(ii) *Let $W_1, \ldots, W_m$ be random variables with the same distribution as $W$. For any $w_n > 0$,*

$$\mathrm{P}\left(w_n > \max_{1 \leq j \leq m}|W_j|\right) \geq 1 - \frac{(\log 2)^{1/2} K(\log m)^{1/2}}{w_n}$$

*for a constant $K$ not depending on $n$.*

PROOF. (i) Without loss of generality, assume $c_i \neq 0, i = 1, \ldots, n$. First, because $\varepsilon_i$ is sub-Gaussian, its Orlicz norm $\|\varepsilon_i\|_{\psi_2} \leq [(1+K)/C]^{1/2}$ [Lemma 2.2.1, Van der Vaart and Wellner (1996)]. By Proposition A.1.6 of Van der Vaart and Wellner (1996), there exists a constant $K_2$ such that

$$\left\|\sum_{i=1}^n c_i \varepsilon_i\right\|_{\psi_2} \leq K_2 \left\{ \mathrm{E}\left|\sum_{i=1}^n c_i \varepsilon_i\right| + \left[\sum_{i=1}^n \|c_i \varepsilon_i\|_{\psi_2}^2\right]^{1/2}\right\}$$

$$\leq K_2 \left\{\sigma + \left[(1+K)C^{-1}\sum_{i=1}^n c_i^2\right]^{1/2}\right\}$$

$$= K_2[\sigma + ((1+K)C^{-1})^{1/2}].$$

(ii) By Lemma 2.2.2 of Van der Vaart and Wellner (1996), $\|\max_{1 \leq j \leq q_n} W_i\|_{\psi_2} \leq K(\log m)^{1/2}$ for a constant $K$. Because $\mathrm{E}|W| \leq (\log 2)^{1/2}\|W\|_{\psi_2}$ for any ran-

dom variable $W$, we have

$$\mathrm{E}\Big(\max_{1\leq j\leq m_n} |W_j|\Big) \leq (\log 2)^{1/2} K(\log m)^{1/2}$$

for a constant $K$. By the Markov inequality, we have

$$\mathrm{P}\Big(w_n > \max_{1\leq j\leq m} |W_j|\Big) = 1 - \mathrm{P}\Big(\max_{1\leq j\leq m_n} |W_j| \geq w_n\Big)$$

$$\geq 1 - \frac{(\log 2)^{1/2} K(\log m)^{1/2}}{w_n}.$$

This completes the proof. $\square$

PROOF OF THEOREM 3. Recall $\xi_{nj} = n^{-1}\sum_{i=1}^{n}(\mathbf{w}_i'\boldsymbol{\beta}_{10})x_{ij}$ as defined in (5). Let $\mathbf{a}_j = (x_{1j},\ldots,x_{nj})'$. Write

$$U_n(\boldsymbol{\beta}) = \sum_{j=1}^{p_n}\sum_{i=1}^{n}(y_i - x_{ij}\beta_j)^2 + \lambda_n\sum_{j=1}^{p_n}|\beta_j|^\gamma$$

$$= \sum_{j=1}^{p_n}\Bigg[\sum_{i=1}^{n}\varepsilon_i^2 + n\beta_j^2 - 2(\boldsymbol{\varepsilon}_n'\mathbf{a}_j + n\xi_{nj})\beta_j + \lambda_n|\beta_j|^\gamma\Bigg].$$

So minimizing $U_n$ is equivalent to minimizing $\sum_{j=1}^{p_n}[n\beta_j^2 - 2(\boldsymbol{\varepsilon}_n'\mathbf{a}_j + n\xi_{nj})\beta_j + \lambda_n|\beta_j|^\gamma]$. Let

$$g_j(\beta_j) \equiv n\beta_j^2 - 2(\boldsymbol{\varepsilon}_n'\mathbf{a}_j + n\xi_{nj})\beta_j + \lambda_n|\beta_j|^\gamma, \qquad j = 1,\ldots,p_n.$$

By Lemma 3, $\beta_j = 0$ is the only solution to $g_j(\beta_j) = 0$ if and only if

$$n^{-1}\lambda_n > c_\gamma(n^{-1}|\boldsymbol{\varepsilon}_n'\mathbf{a}_j + n\xi_{nj}|)^{2-\gamma}.$$

Let $w_n = c_\gamma^{-1/(2-\gamma)}(\lambda_n/n^{\gamma/2})^{1/(2-\gamma)}$. This inequality can be written

(16) $$w_n > n^{-1/2}|\boldsymbol{\varepsilon}_n'\mathbf{a}_j + n\xi_{nj}|.$$

To prove the theorem, it suffices to show that

(17) $$\mathrm{P}\Big(w_n > n^{-1/2}\max_{j\in J_n}|\boldsymbol{\varepsilon}_n'\mathbf{a}_j + n\xi_{nj}|\Big) \to 1$$

and

(18) $$\mathrm{P}\Big(w_n > n^{-1/2}\min_{j\in K_n}|\boldsymbol{\varepsilon}_n'\mathbf{a}_j + n\xi_{nj}|\Big) \to 0.$$

We first prove (17). By condition (B2)(a), there exists a constant $c_0 > 0$ such that

$$\Bigg|n^{-1/2}\sum_{i=1}^{n}x_{ij}x_{ik}\Bigg| \leq c_0, \qquad j \in J_n, k \in K_n,$$

for all $n$ sufficiently large. Therefore,

$$n^{1/2}|\xi_{nj}| = n^{-1/2}\left|\sum_{k=1}^{k_n}\sum_{i=1}^{n} x_{ik}x_{ij}\beta_{0k}\right|$$

(19)
$$\leq n^{-1/2} b_1 \sum_{l=1}^{k_n}\left|\sum_{i=1}^{n} x_{ik}x_{ij}\right|$$

$$\leq b_1 c_0 k_n,$$

where $b_1$ is given in condition (B4). Let $c_1 = b_1 c_0$. By (16) and (19), we have

$$P\left(w_n > n^{-1/2}\max_{j\in J_n}|\varepsilon'_n \mathbf{a}_j + n\xi_{nj}|\right)$$

$$\geq P\left(w_n > n^{-1/2}\max_{j\in J_n}|\varepsilon'_n \mathbf{a}_j| + n^{1/2}\max_{j\in J_n}|\xi_{nj}|\right)$$

(20)
$$\geq P\left(w_n > n^{-1/2}\max_{j\in J_n}|\varepsilon'_n \mathbf{a}_j| + c_1 k_n\right)$$

$$= 1 - P\left(n^{-1/2}\max_{j\in J_n}|\varepsilon'_n \mathbf{a}_j| \geq w_n - c_1 k_n\right)$$

$$\geq 1 - \frac{E(n^{-1/2}\max_{j\in J_n}|\varepsilon'_n \mathbf{a}_j|)}{w_n - c_1 k_n}.$$

By Lemma 4(i), $n^{-1/2}\varepsilon'_n \mathbf{a}_j$ is sub-Gaussian, $1 \leq j \leq m_n$. By condition (B3)(a),

(21)
$$\frac{k_n}{w_n} = \left(\frac{k_n^{(2-\gamma)}}{\lambda_n n^{-\gamma/2}}\right)^{1/(2-\gamma)} \to 0.$$

Thus, by Lemma 4(ii), combining (20) and (21), and by condition (B3)(b),

$$P\left(w_n > n^{-1/2}\max_{j\in J_n}|\varepsilon'_n \mathbf{a}_j + n\xi_{nj}|\right) \geq 1 - \frac{(\log 2)^{1/2} K(\log m_n)^{1/2}}{w_n - c_1 k_n} \to 1.$$

This proves (17). We now prove (18). We have

$$P\left(w_n > \min_{j\in K_n}|n^{-1/2}\varepsilon'_n \mathbf{a}_j + n^{1/2}\xi_{nj}|\right)$$

(22)
$$= P\left(\bigcup_{j\in K_n}\{|n^{-1/2}\varepsilon'_n \mathbf{a}_j + n^{1/2}\xi_{nj}| < w_n\}\right)$$

$$\leq \sum_{j\in K_n} P(|n^{-1/2}\varepsilon'_n \mathbf{a}_j + n^{1/2}\xi_{nj}| < w_n).$$

Write

$$P(|n^{-1/2}\varepsilon'_n \mathbf{a}_j + n^{1/2}\xi_{nj}| < w_n)$$

(23)
$$= 1 - \mathrm{P}(|n^{-1/2}\boldsymbol{\varepsilon}_n'\mathbf{a}_j + n^{1/2}\xi_{nj}| \geq w_n).$$

By condition (B2)(b), $\min_{j \in K_n} |\xi_{nj}| \geq \xi_0 > 0$ for all $n$ sufficiently large. By Lemma 4, $n^{-1/2}\boldsymbol{\varepsilon}_n'\mathbf{a}_j$ are sub-Gaussian. We have

(24)
$$\begin{aligned}
\mathrm{P}(|n^{-1/2}\boldsymbol{\varepsilon}_n'\mathbf{a}_j + n^{1/2}\xi_{nj}| &\geq w_n) \\
&\geq \mathrm{P}(n^{1/2}|\xi_{nj}| - n^{-1/2}|\boldsymbol{\varepsilon}_n'\mathbf{a}_i| \geq w_n) \\
&= 1 - \mathrm{P}(n^{-1/2}|\boldsymbol{\varepsilon}_n'\mathbf{a}_i| > n^{1/2}|\xi_{nj}| - w_n) \\
&\geq 1 - K\exp[-C(n^{1/2}\xi_0 - w_n)^2].
\end{aligned}$$

By (22), (23) and (24), we have

$$\mathrm{P}\bigg(w_n > \min_{j \in K_n} |n^{-1/2}\boldsymbol{\varepsilon}_n'\mathbf{a}_j + n^{1/2}\xi_{nj}|\bigg)$$
$$\leq k_n K \exp[-C(n^{1/2}\xi_0 - w_n)^2].$$

By condition (B3)(a), we have

$$\frac{w_n}{n^{1/2}} = O(1)\bigg(\frac{\lambda_n n^{-\gamma/2}}{n^{(2-\gamma)/2}}\bigg)^{1/(2-\gamma)}$$
$$= O(1)(\lambda_n/n)^{1/(2-\gamma)} = o(1).$$

Therefore,

$$\mathrm{P}\bigg(w_n > \min_{j \in K_n} |n^{-1/2}\boldsymbol{\varepsilon}_n'\mathbf{a}_j + n^{1/2}\xi_{nj}|\bigg) = O(1)k_n \exp(-Cn) = o(1),$$

where the last equality follows from condition (B3)(a). Thus, (18) follows. This completes the proof of Theorem 3. $\square$

PROOF OF THEOREM 4. By Theorem 3, Conditions (B1) to (B4) ensure that the marginal bridge estimator correctly selects covariates with nonzero and zero coefficients with probability converging to one. Therefore, for asymptotic analysis, the second step estimator $\widehat{\boldsymbol{\beta}}_n^*$ can be defined as the value that minimizes $U_n^*$ defined by (6). We now can prove Theorem 4 in two steps. First, under conditions (B1)(a) and (B6), consistency of $\widehat{\boldsymbol{\beta}}_{1n}^*$ follows from the same argument as in the proof of Theorem 1. Then under conditions (B1)(a), (B5) and (B6), asymptotic normality can be proved the same way as in the proof of Theorem 2. This completes the proof of Theorem 4. $\square$

**Acknowledgments.** The authors wish to thank two anonymous referees and the Editor for their helpful comments.

J. Huang  
Department of Statistics  
and Actuarial Science, 241 SH  
University of Iowa  
Iowa City, Iowa 52242  
USA  
E-mail: jian@stat.uiowa.edu

J. L. Horowitz  
Department of Economics  
Northwestern University  
2001 Sheridan Road  
Evanston, Illinois 60208  
USA  
E-mail: joel-horowitz@northwestern.edu

S. Ma  
Division of Biostatistics  
Department of Epidemiology and Public Health  
Yale University  
New Haven, Connecticut 06520  
USA  
E-mail: shuangge.ma@yale.edu